\documentclass[10pt]{amsart}
\usepackage{amssymb,amsmath,amsthm,enumerate,amscd,graphicx,color}
 \usepackage[colorlinks=true, pdfstartview=FitV, linkcolor=blue, citecolor=blue, urlcolor=blue]{hyperref}
\swapnumbers
\newtheorem{thm}{Theorem}[subsection]
\newtheorem{lem}[thm]{Lemma}
\newtheorem{cor}[thm]{Corollary}
\newtheorem{prop}[thm]{Proposition}
\newtheorem{rem}[thm]{Remark}

\newcommand{\thmref}[1]{Theorem~\ref{#1}}
\newcommand{\propref}[1]{Proposition~\ref{#1}}

\newcommand{\lemref}[1]{Lemma~\ref{#1}}
\newcommand{\eqnref}[1]{~(\ref{#1})}

\begin{document}



\title{}

\title{Imaginary Verma modules and Kashiwara algebras for $U_q(\widehat{\mathfrak{sl}(2)})$.}
\author{ Ben Cox}
\author{Vyacheslav Futorny}
\author{Kailash C. Misra}
\keywords{Quantum affine algebras,  Imaginary Verma modules, Kashiwara algebras, simple modules}
\address{Department of Mathematics \\
University of Charleston \\
66 George Street  \\
Charleston SC 29424, USA}\email{coxbl@cofc.edu}
\address{Department of Mathematics\\
 University of S\~ao Paulo\\
 S\~ao Paulo, Brazil}
 \email{futorny@ime.usp.br}
 \address{Department of Mathematics\\
 North Carolina State University\\
 Raleigh, NC 27695-8205, USA}
 \begin{abstract}  We consider imaginary Verma modules for quantum affine algebra
 $U_q(\widehat{\mathfrak{sl}(2)})$  and construct Kashiwara type operators
 and the Kashiwara algebra $\mathcal K_q$. We show that  a certain quotient  $\mathcal N_q^-$ of
 $U_q(\widehat{\mathfrak{sl}(2)})$ is a simple $\mathcal K_q$-module.
\end{abstract}
\date{}
\maketitle

\section{Introduction}
Corresponding to the standard partition of the root system of an affine Lie algebra into
set of positive and negative roots we have a standard Borel subalgebra from which we may
induce the standard Verma modules. However, unlike for finite dimensional semisimple Lie algebras
for an affine Lie algebra there exists other closed partitions of the root system which are not equivalent to the usual partition of the root system under the Weyl group action. Corresponding to such non-standard partitions we have non-standard Borel subalgebras from which one may induce other non-standard Verma-type modules and these typically contain both finite and infinite dimensional weight spaces. The classification of closed subsets of the root system for affine Kac-Moody algebras was obtained by Jakobsen and Kac \cite{JK,MR89m:17032}, and independently by Futorny \cite{MR1078876,MR1175820}. A categorical setting for these modules was introduced in \cite{MR95d:17026}, with certain restrictions, and generalized in \cite{MR97c:17036}. For the algebra $ \widehat{\mathfrak{sl}(2)}$, the only non-standard modules of Verma-type are the {\it imaginary Verma modules} \cite{MR95a:17030}.

Drinfeld \cite{MR802128} and Jimbo \cite{MR797001} independently introduced the quantum group $U_q(\mathfrak g)$ as
$q$-deformations of universal enveloping algebras of a symmetrizable Kac-Moody Lie algebra
$\mathfrak g$. For generic $q$, Lusztig \cite{MR954661} showed that integrable highest weight modules of symmetrizable Kac-Moody algebras can be deformed to those over the corresponding quantum groups in such a way that the dimensions of the weight spaces are invariant under the deformation. Following the framework of \cite{MR954661} and \cite{MR1341758}, {\it quantum imaginary Verma modules}
for the quantum group $U_q(\widehat{\mathfrak{sl}(2)})$ were constructed in \cite {MR97k:17014} and it was shown that these modules are deformations of those over the universal enveloping algebra $U(\widehat{\mathfrak{sl}(2)})$ in such a way that the weight multiplicities, both finite and infinite-dimensional, are preserved.

Kashiwara  (\cite{MR1090425, MR1115118}) from algebraic view point and Lusztig \cite{MR1035415} from geometric view point
introduced global crystal base (equivalently, canonical base) for standard Verma modules
$V_q(\lambda)$ and integrable highest weight modules $L_q(\lambda)$ independently. The crystal base (\cite{MR1090425, MR1115118}) can be thought of as the $q=0$ limit of the global crystal base or canonical base. An important ingredient in the construction of crystal base by Kashiwara in
 \cite{MR1115118}, is a subalgebra $\mathcal {B}_q$ of the quantum group $U_q(\mathfrak g)$ which acts on the negative part $U^-_q(\mathfrak g)$ of the quantum group $U_q(\mathfrak g)$ by left multiplication.
This subalgebra $\mathcal {B}_q$, which we call the Kashiwara algebra, played an important role in the definition of the Kashiwara operators which defines the crystal base.

In this paper we construct an analog of Kashiwara algebra $\mathcal K_q$ for the imaginary Verma module $M_q(\lambda)$ for the quantum group $U_q(\widehat{\mathfrak{sl}(2)})$. Then we prove that certain quotient $\mathcal N_q^-$ of  $U_q(\widehat{\mathfrak{sl}(2)})$ is a simple $\mathcal K_q$-module.
In Sections 2 and 3 we recall necessary definitions and some new results that we need. In Section 4 we define certain operators we call $\Omega$-operators acting on $\mathcal N_q^-$ and prove generalized commutation relations among them. We define the Kashiwara algebra $\mathcal K_q$ in Section 5 and show that  $\mathcal N_q^-$ is a left $\mathcal K_q$-module and define a symmetric invariant bilinear form on $\mathcal N_q^-$. The main result in Section 6 is that for any weight $\lambda$ of level zero the reduced imaginary Verma module $\tilde{M}_q(\lambda)$ is simple if and only if $\lambda(h)\not= 0$ which shows that Lusztig's deformation functor preserves module structure in the case of imaginary Verma modules (see \cite {MR95a:17030}). Finally, in Section 7 we prove that $\mathcal N_q^-$ is
simple as a $\mathcal K_q$-module and that the form defined in Section 5 is nondegenerate.


\section{ Imaginary Verma Modules for $A_1^{(1)}$}

We begin by recalling some basic facts and constructions for the affine
Kac-Moody algebra $ A_1^{(1)}$
and its imaginary Verma modules.
See \cite{K} for Kac-Moody algebra terminology and standard notations.

\subsection{}
Let $\mathbb F$ be a field of characteristic 0. The algebra $A_1^{(1)}$ is the affine Kac-Moody algebra over field $\mathbb F$  with
generalized Cartan matrix $A=(a_{ij})_{0\le i,j \le 1} = \begin{pmatrix}2&-2 \\ -2&\
2 \\ \end{pmatrix}$.
The algebra $A_1^{(1)}$ has a Chevalley-Serre presentation with generators
$e_0, e_1, f_0, f_1, h_0, h_1, d$ and relations
\begin{align*}
&[h_i,h_j] =0, \ \ [h_i,d]=0, \\
&[e_i,f_j] = \delta_{ij}h_i, \\
&[h_i, e_j] = a_{ij}e_j, \ \ [h_i, f_j] = -a_{ij}f_j, \\
&[d, e_j]= \delta_{0,j} e_j, \ \ [d, f_j]=-\delta_{0,j} f_j, \\
&(\text{ad}\, e_i)^3e_j = (\text{ad}\, f_i)^3f_j = 0, \quad i \neq j.
\end{align*}
Alternatively, we may realize $A_1^{(1)}$ through the loop algebra
construction
$$
A_1^{(1)}
\cong {\mathfrak sl}_2 \otimes \mathbb F[t,t^{-1}] \oplus \mathbb F c \oplus
\mathbb F d
$$
with Lie bracket relations
\begin{align*}[x\otimes t^n,y\otimes t^m]&= [x,y] \otimes t^{n+m}
+ n \delta_{n+m,0}(x,y)c, \\
[x \otimes t^n,c] = 0 = [d,c] , \quad & \quad [d,x \otimes t^n] = nx \otimes t^n,
\end{align*}for $x,y \in {\mathfrak sl_2}$, $n,m \in \mathbb Z$, where $(\ , \ )$ denotes the
Killing form on ${\mathfrak sl_2}$.
For $x \in {\mathfrak sl_2}$ and $n \in \mathbb Z$, we write $x(n)$ for
$x \otimes t^n$.

Let $\Delta$ denote the root system of $A_1^{(1)}$, and let
$\{ \alpha_0, \alpha_1\}$ be a basis for $\Delta$.  Let $\delta = \alpha_0 + \alpha_1$,
the minimal imaginary root.  Then
$$
\Delta = \{ \pm \alpha_1 + n\delta\ |\ n \in \mathbb Z\} \cup \{ k\delta\ |\ k \in \mathbb Z
\setminus \{ 0 \} \}.
$$

\subsection{}
The universal enveloping
algebra $
U(A_1^{(1)})$ of $A_1^{(1)}$
is the associative algebra
over $\mathbb F$ with 1
generated by the elements $h_0, h_1, d, e_0, e_1, f_0, f_1$
with defining relations
\begin{align*}&[h_0,h_1]=[h_0,d] = [h_1,d]=0, \\
&h_ie_j-e_jh_i = a_{ij}e_j, \quad h_if_j-f_jh_i=-a_{ij}f_j, \\
&d e_j-e_j d =\delta_{0,j}e_j, \quad d f_j-f_j d = -\delta_{0,j}f_j, \\
&e_if_j-f_je_i = \delta_{ij}h_i, \\
&e_je_i^3-3e_ie_je_i^2+3e_i^2e_je_i-e_i^3e_j = 0 \text{ for } i \neq j, \\
&f_jf_i^3-3f_if_jf_i^2+3f_i^2f_jf_i-f_i^3f_j = 0 \text{ for } i \neq j.
\end{align*}
 Corresponding to the loop algebra formulation of $A_1^{(1)}$ is an
alternative description of
$U(A_1^{(1)})$ as the associative algebra over $\mathbb F$ with 1 generated
by the elements $e(k), f(k)$ $(k \in \mathbb Z)$, $h(l)$ $(l \in
\mathbb Z\setminus \{ 0\})$,
$c, d, h$, with relations
\begin{align*}& [c,u]=0 \ \ \text {for all} \ u\in U(A_1^{(1)}), \\
& [h(k), h(l)]=2k \delta_{k+l,0} c, \\
& [h,d]=0, \ \ [h, h(k)]=0, \\
& [d,h(l)]=l h(l), \ \ [d,e(k)]=ke(k), \ \ [d,f(k)]=kf(k),\\
& [h,e(k)]=2e(k), \ \ [h,f(k)]=-2f(k), \\
& [h(k), e(l)]=2e(k+l), \ \ [h(k), f(l)]=-2f(k+l), \\
& [e(k), f(l)]=h(k+l)+k \delta_{k+l,0} c.
\end{align*}

\subsection{}  A
subset $S$ of the root system $\Delta$ is called {\it closed}
if $\alpha, \beta \in S$ and
$\alpha+\beta \in \Delta$
implies $\alpha+\beta \in S$.  The subset $S$ is called a {\it closed
partition } of the roots if $S$ is closed,
$S \cap(-S) = \emptyset$, and $S\cup -S = \Delta$ \cite{JK},\cite{MR89m:17032},\cite{MR1078876},\cite{MR1175820}.
The set
$$
S= \{ \alpha_1+k\delta \ |\ k\in \mathbb Z \} \cup \{l\delta\ |\ l \in \mathbb Z_{>0} \}
$$
is a closed partition of $\Delta$ and is $W\times
\{\pm{1}\}$-inequivalent to the standard
partition of the root system into positive and negative roots \cite{MR95a:17030}.

For ${\mathfrak g}= A_1^{(1)}$, let ${\mathfrak g}_{\pm}^{(S)}=\sum_{\alpha \in S}
{\mathfrak g}_{\pm \alpha}$.  In the loop algebra formulation of ${\mathfrak g}$, we have
that ${\mathfrak g}_+^{(S)}$ is the
subalgebra generated by $e(k)$ $(k \in \mathbb Z)$ and $h(l)$ $(l\in \mathbb Z_{>0})$
and ${\mathfrak g}_-^{(S)}$ is the subalgebra generated by $f(k)$ $(k \in \mathbb Z)$ and
$h(-l)$ $(l\in \mathbb Z_{>0})$.  Since $S$ is a partition of the root system,
the algebra has a direct
 sum decomposition
$$
{\mathfrak g}={\mathfrak g}_{-}^{(S)} \oplus {\mathfrak h} \oplus {\mathfrak g}_{+}^{(S)}.
$$
Let $U({\mathfrak g}_{\pm}^{(S)})$ be the universal enveloping algebra of
${\mathfrak g}_{\pm}^{(S)}$. Then, by the PBW theorem, we have
$$
U({\mathfrak g}) \cong U({\mathfrak g}_{-}^{(S)}) \otimes U({\mathfrak h})\otimes U({\mathfrak g}_{+}^{(S)}),
$$
where $U({\mathfrak g}_{+}^{(S)})$ is generated by $ e(k)$ $(k\in \mathbb Z)$, $h(l)$
$(l\in \mathbb Z_{>0})$,
$U({\mathfrak g}_{-}^{(S)})$ is generated by $f(k)$ $(k\in \mathbb Z)$, $h(-l)$ $(l\in
\mathbb Z_{>0})$ and $U({\mathfrak h})$,
the universal enveloping algebra of ${\mathfrak h}$, is generated by
$h,c$ and $d$.

Let $\lambda\in P$, the weight lattice of ${\mathfrak g}=A_1^{(1)}$.
A $U({\mathfrak g})$-module $V$ is called a {\it weight} module if
$V=\oplus_{\mu \in P} V_{\mu}$, where
$$
V_{\mu}=\{ v \in V\ |\ h\cdot v=\mu(h)v, c\cdot v=\mu(c)v,
d\cdot v = \mu(d)v \}.
$$
Any submodule of a weight module is a weight module.
A $U({\mathfrak g})$-module $V$ is
called an {\it $S$-highest weight module}
with highest weight $\lambda$ if there is a
non-zero $v_{\lambda} \in V$ such that
(i) $u^+ \cdot v_{\lambda} = 0$ for all $u^+\in U({\mathfrak g}_{+}^{(S)})
\setminus \mathbb F^*$, (ii) $h\cdot v_{\lambda}=\lambda(h)v_{\lambda}$, $c\cdot v_{\lambda}
= \lambda(c)v_{\lambda}$, $d\cdot v_{\lambda}  = \lambda(d)v_{\lambda}$,
(iii) $V=U({\mathfrak g})\cdot v_{\lambda} = U({\mathfrak g}_{-}^{(S)}) \cdot v_{\lambda}$.
An $S$-highest weight module is a weight module.

For $\lambda \in P$, let $I_S(\lambda)$ denote the ideal of $U(A_1^{(1)})$
generated by
$e(k)$ $(k\in \mathbb Z)$, $h(l)$ $(l>0)$, $h-\lambda(h) 1$,
$c-\lambda(c) 1$, $d-\lambda(d) 1$.
Then we define $M(\lambda) = U(A_1^{(1)})/I_S(\lambda)$ to be the {\it
imaginary Verma module} of
$A_1^{(1)}$ with highest weight $\lambda$.
Imaginary Verma modules have many structural features similar to those of
standard
Verma modules,
with the exception of the infinite-dimensional weight spaces.
Their properties were investigated in \cite{MR95a:17030}, from which we recall
 the following
proposition \cite[ Proposition 1, Theorem 1]{MR95a:17030}.

\begin{prop} (i) $M(\lambda)$ is a $U({\mathfrak g}_-^{(S)})$-free
module of rank 1 generated by
the $S$-highest weight vector $1\otimes 1$ of weight $\lambda$.\newline
(ii) $\dim M(\lambda)_{\lambda} =1$; $0<\dim M(\lambda)_{\lambda - k \delta} < \infty$
for any integer $k>0$; if $\mu \neq \lambda- k\delta$ for any integer
 $k \ge 0$ and
$M(\lambda)_{\mu} \neq 0$, then $\dim M(\lambda)_{\mu} = \infty$.\newline
(iii) Let $V$ be a $U(A_1^{(1)})$-module generated by some $S$-highest
weight vector $v$
 of weight $\lambda$.  Then there exists a unique surjective homomorphism
$\varphi:M(\lambda) \to V$ such that $\varphi(1 \otimes 1) =v$. \newline
(iv) $M(\lambda)$ has a unique maximal submodule. \newline
(v) Let $\lambda, \mu \in P$.  Any non-zero element of
$\text{Hom}\,_{U(A_1^{(1)})}(M(\lambda),M(\mu))$ is
injective.\newline
(vi) $M(\lambda)$ is irreducible if and only if $\lambda(c)\neq 0$.
\hskip 1cm $\square$
\end{prop}

\section{The quantum group $U_q(A_1^{(1)})$}

\subsection{}  The {\it quantum group}
$U_q(A_1^{(1)})$ is the $\mathbb F(q^{1/2})$-algebra with 1 generated by
$$ e_0, e_1, f_0, f_1, K_0^{\pm 1}, K_1^{\pm 1}, D^{\pm 1} $$
with defining relations:
\begin{align*}& DD^{-1}=D^{-1}D=K_iK_i^{-1}=K_i^{-1}K_i=1, \\
& e_if_j-f_je_i = \delta_{ij}\frac{K_i-K_i^{-1}}{q-q^{-1}}, \\
& K_ie_iK_i^{-1}=q^2e_i, \ \ K_i f_i K_i^{-1} =q^{-2}f_i, \\
& K_i e_jK_i^{-1} = q^{-2}e_j, \ \
K_i f_jK_i^{-1} = q^2f_j, \quad i\neq j, \\
& K_iK_j-K_jK_i = 0, \ \ K_iD-DK_i=0, \\
& De_iD^{-1}=q^{\delta_{i,0}} e_i, \ \
Df_iD^{-1}=q^{-\delta_{i,0}} f_i, \\
& e_i^3e_j-[3]e_i^2e_je_i+[3]e_ie_je_i^2-e_je_i^3 =0, \quad i\neq j, \\
& f_i^3f_j-[3]f_i^2f_jf_i+[3]f_if_jf_i^2-f_jf_i^3 = 0, \quad i\neq j, \\
\end{align*}where, $[n] = \frac{q^n-q^{-n}}{q-q^{-1}}$.

The quantum group $U_q(A_1^{(1)})$
can be given a Hopf algebra structure with a comultiplication given by
\begin{align*}
\Delta(K_i) &= K_i \otimes K_i, \\
\Delta(D)&=D\otimes D, \\
\Delta(e_i) &= e_i\otimes K_i^{-1} + 1\otimes e_i, \\
\Delta(f_i) &= f_i\otimes 1 + K_i \otimes f_i, \\
\end{align*}and an antipode given by
\begin{align*}
s(e_i) &=-e_iK_i^{-1}, \\
s(f_i) &= -K_if_i, \\
s(K_i) &= K_i^{-1},\\
s(D) &= D^{-1}.
\end{align*}

  There is an alternative realization for $U_q(A_1^{(1)})$,
due to Drinfeld
\cite{MR802128}, which we shall also need.  Let
$U_q$ be the associative algebra with $1$ over $\mathbb F(q^{1/2})$
generated by the
elements $x^{\pm }(k)$ ($k\in \mathbb Z$), $a(l)$ ($l \in \mathbb Z
\setminus \{0\}$), $K^{\pm 1}$,
$D^{\pm 1}$, and $\gamma^{\pm \frac12}$ with the following defining
relations:
\begin{align}
DD^{-1}=D^{-1}D&=KK^{-1}=K^{-1}K=1,  \\
[\gamma^{\pm \frac 12},u] &= 0 \quad \forall u \in U, \\
[a(k),a(l)] &= \delta_{k+l,0} \frac{[2k]}{k} \frac{\gamma^k -
\gamma^{-k}}{q-q^{-1}},  \\
[a(k),K]&=0,\quad [D,K]=0,  \\
D a(k)D^{-1}&=q^k a(k),  \\
Dx^{\pm}(k)D^{-1}&=q^{ k}x^{\pm}(k),  \\
Kx^{\pm}(k)K^{-1} &= q^{\pm 2}x^{\pm}(k),   \\
[a(k),x^{\pm}(l)]&= \pm \frac{[2k]}{k}\gamma^{\mp \frac{|k|}{2}}x^{\pm}(k+l), \label{axcommutator}  \\
    x^{\pm}(k+1)x^{\pm}(l) &- q^{\pm 2}
x^{\pm}(l)x^{\pm}(k+1)\label{Serre}   \\
&= q^{\pm 2}x^{\pm}(k)x^{\pm}(l+1)
    - x^{\pm}(l+1)x^{\pm}(k),\notag \\
[x^+(k),x^-(l)]&=
    \frac{1}{q-q^{-1}}\left( \gamma^{\frac{k-l}{2}}\psi(k+l) -
    \gamma^{\frac{l-k}{2}}\phi(k+l)\right), \label{xcommutator}   \\
\text{where  }
\sum_{k=0}^{\infty}\psi(k)z^{-k} &= K \exp\left(
(q-q^{-1})\sum_{k=1}^{\infty} a(k)z^{-k}\right),\\
\sum_{k=0}^{\infty}
\phi(-k)z^k&= K^{-1} \exp\left( - (q-q^{-1})\sum_{k=1}^{\infty}
a(-k)z^k\right).
\end{align}
The algebras $U_q(A_1^{(1)})$
and $U_q$ are isomorphic \cite{MR802128}. The action of the
isomorphism,
which we shall call the {\it Drinfeld Isomorphism,} on the generators of
$U_q(A_1^{(1)})$ is:
\begin{align*}e_0 &\mapsto x^-(1)K^{-1}, \ \ f_0 \mapsto Kx^+(-1), \\
e_1 &\mapsto x^+(0), \ \ f_1 \mapsto x^-(0), \\
K_0 &\mapsto \gamma K^{-1}, \ \ K_1 \mapsto K, \ \ D \mapsto D.
\end{align*}
If one uses the formal sums
\begin{equation}
\phi(u)=\sum_{p\in\mathbb Z} \phi(p)u^{-p},\enspace \psi(u)=\sum_{p\in\mathbb Z}\psi(p)u^{-p},\enspace
x^{\pm }(u)=\sum_{p\in\mathbb Z} x^\pm(p)u^{-p}
\end{equation}
Drinfeld's relations (3), (8)-(10) can be written as
\begin{gather}
[\phi(u),\phi(v)]=0=[\psi(u),\psi(v)] \\
\phi(u)x^\pm (v)\phi(u)^{-1}=g(uv^{-1}\gamma^{\mp 1/2})^{\pm 1}x^\pm (v)\label{phix} \\
\psi(u)x^\pm (v)\psi(u)^{-1}=g(vu^{-1}\gamma^{\mp 1/2})^{\mp 1}x^\pm (v)\label{psix} \\
(u-q^{\pm2} v)x^\pm (u)x^\pm (v)=(q^{\pm 2}u-v)x^\pm(v)x^\pm(u) \\
[x^+(u),x^-(v)]=(q-q^{-1})^{-1}(\delta(u/v\gamma)\psi(v\gamma^{1/2})-\delta(u\gamma/v)\phi(u\gamma^{1/2}))\label{xx}
\end{gather}
where $g(t)=g_q(t)$ is the Taylor series at $t=0$ of the function $(q^2t-1)/(t-q^2)$ and $\delta(z)=\sum_{k\in\mathbb Z}z^{k}$ is the formal Dirac delta function.
\begin{rem}
Writing $g(t)=g_q(t)=\sum_{p\geq 0}g(p)t^p$ we have
$$g(0)=q^{-2}, \quad g(p)=(1-q^4)q^{-2p-2}, \quad p>0.$$
Note that $g_q(t)^{-1}=g_{q^{-1}}(t).$
\end{rem}

We will need the following identity later:
\begin{lem}
\begin{equation}
\label{quteformula}
\text{\rm exp}\,\left((q-q^{-1})
\sum_{k=1}^\infty \frac{-[2k]}{k}z^{-k}\right)=1+(1-q^4)\sum_{r=1}^\infty \left(zq^2\right)^{-r}=q^{2}g(1/z)
\end{equation}
\end{lem}
\begin{proof}
\begin{align*}
\text{\rm exp}\,\left((q-q^{-1})
\sum_{k=1}^\infty \frac{-[2k]}{k}z^{-k}\right)&=\text{\rm exp}\,\left(\sum_{k=1}^\infty \frac{1}{k}(zq^2)^{-k}-\sum_{k=1}^\infty \frac{1}{k}\left(\frac{z}{q^2}\right)^{-k}\right)  \\
&=\text{\rm exp}\,\left(\sum_{k=1}^\infty \frac{1}{k}(zq^2)^{-k}\right)\text{\rm exp}\,\left(-\sum_{k=1}^\infty \frac{1}{k}\left(\frac{z}{q^2}\right)^{-k}\right)  \\
&=\left(\frac{1}{1-\frac{1}{zq^2}}\right)\left(1-\frac{q^2}{z}\right)\\
&=\sum_{k=0}^\infty \left(\frac{1}{zq^2}\right)^k-q^4\sum_{k=1}^\infty \left(\frac{1}{zq^2}\right)^{k}
\end{align*}
\end{proof}

\subsection{}
Using the root partition $ S = \{\alpha_1 + k\delta  \ |\ k \in \mathbb Z \} \cup
\{ l\delta \ |\ l \in \mathbb Z_{>0}\}$ from Section 2.3, we define:
\newline
$U_q^+(S)$ to be the subalgebra of $U_q$ generated by $x^+(k)$ $(k\in \mathbb Z)$ and
$a(l)$ $(l>0)$;
\newline
$U_q^-(S)$ to be the subalgebra of $U_q$ generated by $x^-(k)$
$(k \in \mathbb Z )$ and
$a(-l)$ $(l>0)$, and
\newline
$U_q^0(S)$ to be the subalgebra of $U_q$ generated by $K^{\pm 1}$,
$\gamma^{\pm 1/2}$, and $D^{\pm 1}$.

Then we have the following PBW theorem.
\begin{thm}[\cite{MR97k:17014}] \label{PBW}A basis for $U_q$ is the set of monomials
of the form
$$
x^-a^-K^\alpha D^{\beta}\gamma^{\mu/2}a^+x^+
$$
where
\begin{align*}
x^\pm&=x^\pm(m_1)^{n_1}\cdots x^\pm(m_k)^{n_k},
        \quad\quad m_i<m_{i+1},\quad m_i\in\mathbb Z,  \\
a^\pm&=a(r_1)^{s_1}\cdots a(r_l)^{s_l},
         \quad\quad r_i<r_{i+1},\quad \pm r_i\in\mathbb N^*,
\end{align*}and $\alpha,\beta,\mu\in\mathbb Z,\, n_i,s_i\in\mathbb N$.
In particular, $U_q \cong U_q^-(S) \otimes U_q^0(S) \otimes U_q^+(S)$.
\end{thm}

We may order monomials in $u$ in such a way that $-r_1\leq -r_2\leq
\ldots \leq -r_l$ and we compare elements lexicographically.

Considering Serre's relation \eqnref{Serre} with $k=l$, we get
$$
x^{-}(k+1)x^{-}(k) = q^{- 2}x^{-}(k)x^{-}(k+1).
$$
The product on the right side is in the correct order for a basis  element.
If $k+1>l$ and $k\neq l$ in \eqnref{Serre}, then $k+1>  l+1$ so that $k\geq l+1$, and thus we can write
$$
x^{-}(k+1)x^{-}(l) =q^{- 2}x^{-}(l)x^{-}(k+1)  + q^{- 2}x^{-}(k)x^{-}(l+1) - x^{-}(l+1)x^{-}(k)
$$
and then after repeating the above identity (for example the next step is to replace $k+1$ by $k$ and $l$ by $l+1$ on the left), we will eventually arrive at terms that are in the correct order.    In particular if $k+1>l$ and $k\neq l$ note that $x^{-}(l)x^{-}(k+1)<x^{-}(l+1)x^{-}(k)$.

\section{$\Omega$-operators and their relations}

Let $\mathbb N^{\mathbb N^*}$ denote the set of all functions from $\{k\delta\,|\, k\in\mathbb N^*\}$ to $\mathbb N$ with finite support.   Then we can write
$$
a^+=a^{(s_k)}_+:=a (r_1)^{s_1}\cdots a (r_l)^{s_l},\quad a^-:=a^{(s_k)}_-=
a(-r_1)^{s_1}\cdots a(-r_l)^{s_l}
$$
for $f=(s_k)\in\mathbb N^{\mathbb N^*}$ whereby $f(r_k)=s_k$ and $f(t)=0$ for $t\neq r_i, 1\leq i\leq  l$.

Consider now the subalgebra $ \mathcal N_q^-$, generated by $\gamma^{\pm1/2}$, and $x^-(l)$, $l\in\mathbb Z$. Note that the corresponding relations (9) hold in $ \mathcal N_q^-$.

\begin{lem} \label{Plemma} Fix $k\in\mathbb Z$.  Then for any $P\in \mathcal N_q^-$, there exists unique
$$
Q(a,(q_k)),R(c,(r_l))\in \mathcal N_q^-,\quad a,b\in \mathbb Z , (q_l),(r_m)\in\mathbb N^{\mathbb N^*},
$$
  such that
\begin{equation}
[x^-(k),P]=\sum \frac{ a_+^{(q_l)}K^a Q(a,(q_l))}{q-q^{-1}}+\sum \frac{a_-^{(r_m)}K^b\ R(b,(r_m))}{q-q^{-1}}.
\end{equation}

\end{lem}
\begin{proof}  The uniqueness follows from \thmref{PBW} above.
Now any element in $\mathcal N_q^-$ is a sum of products of elements of the form
$$
P=\gamma^{l/2}x^- (m_1) \cdots x^- (m_k) ,\quad \text{ where } m_i\in\mathbb Z, m_1\leq m_2\leq\cdots\leq m_k, k\geq 0,\, l\in\mathbb Z
$$
and such a product is a summand of
$$
P=P(v_1,\dots,v_k):=\gamma^{l/2}x^-(v_1)\cdots x^-(v_k)
$$
Set $\bar P=x^-(v_1)\cdots x^-(v_k)$ and $\bar P_l=x^-(v_{1})\cdots x^-(v_{l-1})x^-(v_{l+1})\cdots x^-(v_k)$.

Then we have by \eqnref{phix} and \eqnref{psix},
\begin{align*}
x^-(v_1)\cdots x^-(v_{l-1})\psi(v_l\gamma^{1/2}) =\prod_{j=1}^{l-1}g(v_jv_l^{-1} )^{-1}
\psi(v_l\gamma^{1/2})x^-(v_{1})\cdots x^-(v_{l-1}) \\
x^-(v_1)\cdots x^-(v_{l-1})\phi(u\gamma^{1/2})=\prod_{j=1}^{l-1}
g(u\gamma v_j^{-1})\phi(u\gamma^{1/2})x^-(v_{1})\cdots x^-(v_{l-1}), \\
\end{align*}
so that by \eqnref{xx}
\begin{align*}
[x^+(u),&x^-(v_1)\cdots x^-(v_k)]=\sum_{l=1}^kx^-(v_1)\cdots [x^+(u),x^-(v_l)]\cdots x^-(v_k)  \\
&=\sum_{l=1}^kx^-(v_1)\cdots\left (\frac{\delta(u/v_l\gamma)\psi(v_l\gamma^{1/2})-\delta(u\gamma/v_l)\phi(u\gamma^{1/2})}{q-q^{-1}}\right)\cdots x^-(v_k)  \\
&=\sum_{l=1}^kx^-(v_1)\cdots x^-(v_{l-1})\psi(v_l\gamma^{1/2})x^-(v_{l+1})\cdots x^-(v_k) \frac{\delta(u/v_l\gamma)}{q-q^{-1}} \\
&\quad-\sum_{l=1}^kx^-(v_1)\cdots x^-(v_{l-1})\phi(u\gamma^{1/2})x^-(v_{l+1})\cdots x^-(v_k) \frac{\delta(u\gamma/v_l)}{q-q^{-1}}\\  \\
&=\sum_{l=1}^k\prod_{j=1}^{l-1}g(v_jv_l^{-1} )^{-1}
\frac{\psi(v_l\gamma^{1/2})\delta(u/v_l\gamma)}{q-q^{-1}}  \bar P_l
\\
&\quad-\sum_{l=1}^k\prod_{j=1}^{l-1}
g(u\gamma v_j^{-1})\frac{\phi(u\gamma^{1/2})\delta(u\gamma/v_l)}{q-q^{-1}}  \bar P_l \\  \\
&=\frac{\psi(u\gamma^{-1/2})}{q-q^{-1}} \sum_{l=1}^k\prod_{j=1}^{l-1}g_{q^{-1}}(v_j/v_l)
  \bar P_l \delta(u/v_l\gamma)    \\
&\quad-\frac{\phi(u\gamma^{1/2})}{q-q^{-1}}
\sum_{l=1}^k\prod_{j=1}^{l-1}
g(v_l/v_j)\bar P_l\delta(u\gamma/v_l) \\
\end{align*}
\end{proof}

\lemref{Plemma} motivates the definition of a family of operators as follows.
Set
$$
G_l=G_l^{1/q}:=\prod_{j=1}^{l-1}g_{q^{-1}}(v_j/v_l),\quad G_l^{q}=\prod_{j=1}^{l-1}g(v_l /v_j)
$$
where $G_1:=1$.
Now define a collection of operators $\Omega_\psi(k),\Omega_\phi(k):\mathcal N_q^-\to \mathcal N_q^-$, $k\in\mathbb Z$, in terms of the generating functions
$$
\Omega_\psi(u)=\sum_{l\in\mathbb Z}\Omega_\psi(l)u^{-l},\quad \Omega_\phi(u)=\sum_{l\in\mathbb Z}\Omega_\phi(l)u^{-l}
$$
by
\begin{align}\label{definingomegapsi}
\Omega_\psi(u)(\bar P):&=\gamma^{m} \sum_{l=1}^kG_l
  \bar P_l \delta(u/v_l\gamma) \\
   \Omega_\phi(u)(\bar P):&=\gamma^{m}\sum_{l=1}^kG_l^q  \bar P_l\delta(u\gamma/v_l).\label{definingomegaphi}
\end{align}
Then we can write the above computation in the proof of Lemma 4.0.2 as
\begin{equation}\label{xplusP}
[x^+(u),\bar P]=(q-q^{-1})^{-1}\left(\psi(u\gamma^{-1/2})\Omega_\psi(u)(\bar P)-\phi(u\gamma^{1/2})\Omega_\phi(u)(\bar P)\right),
\end{equation}

Note that $\Omega_{\psi}(u)(1)=\Omega_{\phi}(u)(1)=0$.  More generally let us write
$$
\bar P=x^-(v_1)\cdots x^-(v_k)=\sum_{n\in\mathbb Z}\sum_{n_1,n_2,\dots,n_k\in\mathbb Z\atop n_1+\cdots +n_k=n}x^-(n_1)\cdots x^-(n_k)v_1^{-n_1}\cdots v_k^{-n_k}
$$
Then
\begin{align*}
\psi&(u\gamma^{-1/2})\Omega_\psi(u)(\bar P)\\
&=\sum_{k\geq 0}\sum_{p\in\mathbb Z}\sum_{n_i\in\mathbb Z } \gamma^{k/2}\psi(k)\Omega_\psi(p)(x^-(n_1)\cdots x^-(n_k))v_1^{-n_1}\cdots v_k^{-n_k}u^{-k-p}  \\
&=\sum_{n_i\in\mathbb Z } \sum_{m\in\mathbb Z}\sum_{k\geq 0}\gamma^{k/2}\psi(k)\Omega_\psi(m-k)(x^-(n_1)\cdots x^-(n_k))v_1^{-n_1}\cdots v_k^{-n_k}u^{-m}
\end{align*}
while
\begin{equation*}
[x^+(u),\bar P]=\sum_{m\in\mathbb Z}\sum_{n_1,n_2,\dots,n_k\in\mathbb Z } [x^+(m),x^-(n_1)\cdots x^-(n_k)]v_1^{-n_1}\cdots v_k^{-n_k}u^{-m}.
\end{equation*}
Thus for a fixed $m$ and $k$-tuple $(n_1,\dots,n_k)$ the sum
$$
\sum_{k\geq 0}\gamma^{k/2}\psi(k)\Omega_\psi(m-k)(x^-(n_1)\cdots x^-(n_k))
$$
must be finite.  Hence
\begin{equation}\label{omegalocalfin}
\Omega_\psi(m-k)(x^-(n_1)\cdots x^-(n_k))=0,
\end{equation}
 for $k$ sufficiently large.

\begin{prop} \label{commutatorprop} Consider $x^-(v)=\sum_mx^-(m)v^{-m}$ as a formal power series of left multiplication operators $x^-(m):\mathcal N_q^-\to \mathcal N_q^-$.  Then
\begin{align}
\Omega_\psi(u)x^-(v)&=\delta(v\gamma/u)+g_{q^{-1}}(v\gamma/u)x^-(v)\Omega_\psi(u),
\label{omegapsi}\\
  \Omega_\phi(u)x^-(v)&=\delta(u\gamma/v)+g(u\gamma/v)x^-(v)\Omega_\phi(u)\label{omegaphi}  \\
(q^2u_1-u_2)\Omega_\psi(u_1)\Omega_\psi(u_2)&=(u_1-q^2u_2)\Omega_\psi(u_2)\Omega_\psi(u_1) \label{psipsi} \\
(q^2u_1-u_2)\Omega_\phi(u_1)\Omega_\phi(u_2)&=(u_1-q^2u_2)\Omega_\phi(u_2)\Omega_\phi(u_1)   \label{phiphi} \\
(q^2\gamma^2u_1-u_2)\Omega_\phi(u_1)\Omega_\psi(u_2)&=(\gamma^2u_1-q^2u_2)\Omega_\psi(u_2)\Omega_\phi(u_1)\label{omegaphipsi}
\end{align}
\end{prop}

\begin{proof} Setting $\bar P=x^-(v_1)\cdots x^-(v_k)$ we get
\begin{align*}
\Omega_\psi(u)x^-(v)&(\bar P)=
  x^-(v_{1}) \cdots x^-(v_k) \delta(u/v\gamma) \\
  &\quad +x^-(v)\sum_{l=1}^k g_{q^{-1}}(v/v_l)G_l
   \bar P_l \delta(u/v_l\gamma)  \\
  &=
 \bar P \delta(u/v\gamma) +x^-(v)g_{q^{-1}}(v\gamma /u)\Omega_\psi(u)\bar P.
\end{align*}
Similarly
\begin{align*}
\Omega_\phi(u)x^-(v)(\bar P)&
    =
 x^-(v_{1}) \cdots x^-(v_k) \delta(u\gamma/v ) \\
      &\quad +x^-(v)\sum_{l=1}^k g (v_l/v)G_l^{q}
     \bar P_l \delta(u\gamma/v_l )  \\
  &=\bar P \delta(v/u\gamma) +x^-(v)g (u\gamma/v)\Omega_\phi(u)\bar P.
\end{align*}
One can prove \eqnref{psipsi} and \eqnref{phiphi} directly from their definitions, \eqnref{definingomegapsi} and \eqnref{definingomegaphi}, but there is another way to prove this identity and it goes as follows:
\begin{align*}
\Omega_\psi(u_1)\Omega_\psi(u_2)&x^-(v) =
  \Omega_\psi(u_1) \delta(v\gamma/u_2) +\Omega_\psi(u_1) x^-(v)g_{q^{-1}}(v\gamma /u_2)\Omega_\psi(u_2)
  \\
  &= \Omega_\psi(u_1) \delta(v\gamma/u_2) +g_{q^{-1}}(v\gamma /u_2)\Omega_\psi (u_2)\delta(v\gamma/u_1)
 \\
  &\quad +g_{q^{-1}}(v\gamma /u_2)g_{q^{-1}}(v \gamma/u_1)x^-(v)\Omega_\psi(u_1)\Omega_\psi (u_2)
\end{align*}
and on the other hand
\begin{align*}
\Omega_\psi(u_2)\Omega_\psi(u_1)&x^-(v) =
  \Omega_\psi(u_2) \delta(v\gamma/u_1) +\Omega_\psi(u_2) x^-(v)g_{q^{-1}}(v\gamma /u_1)\Omega_\psi(u_1)
  \\
  &= \Omega_\psi(u_2) \delta(v\gamma/u_1) +g_{q^{-1}}(v\gamma /u_1)\Omega_\psi (u_1)\delta(v \gamma/u_2)
 \\
  &\quad +g_{q^{-1}}(v\gamma /u_1)g_{q^{-1}}(v \gamma/u_2)x^-(v)\Omega_\psi(u_2)\Omega_\psi (u_1)
\end{align*}

Thus setting $S=(u_1-q^{-2}u_2)\Omega_\psi(u_1)\Omega_\psi(u_2)-(q^{-2}u_1-u_2)\Omega_\psi(u_2)\Omega_\psi (u_1)$ we get
\begin{align*}
Sx^-(v)&=(u_1-q^{-2}u_2) \Omega_\psi(u_1) \delta(v\gamma /u_2) +(u_1-q^{-2}u_2)g_{q^{-1}}(v\gamma /u_2)\Omega_\psi (u_2)\delta(v \gamma/u_1)
 \\
  &\quad +(u_1-q^{-2}u_2)g_{q^{-1}}(v\gamma /u_2)g_{q^{-1}}(v \gamma/u_1)x^-(v)\Omega_\psi(u_1)\Omega_\psi (u_2)  \\
&\quad-(q^{-2}u_1-u_2)\Omega_\psi(u_2) \delta(v\gamma/u_1) -(q^{-2}u_1-u_2)g_{q^{-1}}(v\gamma /u_1)\Omega_\psi (u_1)\delta(v \gamma/u_2)
 \\
  &\quad -(q^{-2}u_1-u_2)g_{q^{-1}}(v\gamma /u_1)g_{q^{-1}}(v \gamma/u_2)x^-(v)\Omega_\psi(u_2)\Omega_\psi (u_1) \\ \\
 &=\left((u_1-q^{-2}u_2)-(q^{-2}u_1-u_2)g_{q^{-1}}(v\gamma /u_1)\right)\Omega_\psi (u_1)\delta(v \gamma/u_2)  \\
&\quad+\left((u_1-q^{-2}u_2)g_{q^{-1}}(v\gamma /u_2) -(q^{-2}u_1-u_2)\right)\Omega_\psi(u_2) \delta(v\gamma/u_1)
 \\
  &\quad+g_{q^{-1}}(v\gamma /u_2)g_{q^{-1}}(v \gamma/u_1)x^-(v) \\
  	&\hskip 50pt \times \left((u_1-q^{-2}u_2)\Omega_\psi(u_1)\Omega_\psi (u_2) -(q^{-2}u_1-u_2) ) \Omega_\psi(u_2)\Omega_\psi (u_1)\right) \\ \\
&=g_{q^{-1}}(v\gamma /u_2)g_{q^{-1}}(v \gamma/u_1)x^-(v)S
\end{align*}
Hence
$$
Sx^-(v_1)\cdots x^-(v_n)=  \prod_{i=1}^ng_{q^{-1}}(v_i \gamma/u_1)g_{q^{-1}}(v_i\gamma /u_2)x^-(v_1)\cdots x^-(v_n) S,
$$
which implies, after applying this to $1$  that  $S=0$.

Next we have
\begin{align*}
\Omega_\phi(u_1)\Omega_\phi(u_2)&x^-(v) =
  \Omega_\phi(u_1) \delta(v/u_2 \gamma) +\Omega_\phi(u_1) x^-(v)g(u_2 \gamma/v)\Omega_\phi(u_2) \\
  &= \Omega_\phi(u_1) \delta(v/u_2 \gamma) +g(u_2 \gamma/v)\Omega_\phi (u_2)\delta(v/u_1 \gamma) \\
  &\quad +g(u_2 \gamma/v)g(u_1 \gamma/v)x^-(v)\Omega_\phi(u_1)\Omega_\phi (u_2)
\end{align*}
and on the other hand
\begin{align*}
\Omega_\phi(u_2)\Omega_\phi(u_1)&x^-(v) =
\Omega_\phi(u_2) \delta(v/u_1 \gamma) +\Omega_\phi(u_2) x^-(v)g(u_1 \gamma/v)\Omega_\phi(u_1) \\
&= \Omega_\phi(u_2) \delta(v/u_1 \gamma) +g(u_1 \gamma/v)\Omega_\phi (u_1)\delta(v/u_2 \gamma)\\
&\quad +g(u_1 \gamma/v)g(u_2 \gamma/v)x^-(v)\Omega_\phi(u_2)\Omega_\phi (u_1)
\end{align*}

So if we set   $S=(u_1-q^{-2}u_2)\Omega_\phi(u_1)\Omega_\phi(u_2)-(q^{-2}u_1-u_2)\Omega_\phi(u_2)\Omega_\phi (u_1)$ we get
\begin{align*}
Sx^-(v)  &= (u_1-q^{-2}u_2)\Omega_\phi(u_1) \delta(v /u_2 \gamma)
    +(u_1-q^{-2}u_2)g(u_2 \gamma/v)\Omega_\phi (u_2)\delta(v/u_1 \gamma)  \\
&\quad +(u_1-q^{-2}u_2)g(u_2 \gamma/v)g(u_1 \gamma/v)x^-(v)\Omega_\phi(u_1)\Omega_\phi (u_2)  \\
&\quad -(q^{-2}u_1-u_2)\Omega_\phi(u_2) \delta(v/u_1 \gamma) -(q^{-2}u_1-u_2)g(u_1 \gamma/v)\Omega_\phi (u_1)\delta(v/u_2 \gamma)\\
&\quad -(q^{-2}u_1-u_2)g(u_1 \gamma/v)g(u_2 \gamma/v)x^-(v)\Omega_\phi(u_2)\Omega_\phi (u_1)
 \\ \\
&= \left((u_1-q^{-2}u_2) -(q^{-2}u_1-u_2)g(u_1 \gamma/v)\right)\Omega_\phi (u_1)\delta(v/u_2 \gamma)      \\
&\quad+\left((u_1-q^{-2}u_2)g(u_2 \gamma/v) -(q^{-2}u_1-u_2)\right) \Omega_\phi (u_2)\delta(v/u_1 \gamma)  \\
&\quad +g(u_2 \gamma/v)g(u_1 \gamma/v)x^-(v) \\
&\hskip 50pt \times\left((u_1-q^{-2}u_2)\Omega_\phi(u_1)\Omega_\phi (u_2) -(q^{-2}u_1-u_2) \Omega_\phi(u_2)\Omega_\phi (u_1)\right)
 \\ \\
&=   g(u_2 \gamma/v)g(u_1 \gamma/v)x^-(v)S.
\end{align*}
As in the calculation for \eqnref{psipsi} we get $S=0$.

Moreover
\begin{align*}
\Omega_\phi(u_1)\Omega_\psi(u_2)&x^-(v) =
  \Omega_\phi(u_1) \delta(v\gamma/u_2) +\Omega_\phi(u_1) x^-(v)g_{q^{-1}}(v\gamma /u_2)\Omega_\psi(u_2)
  \\
  &= \Omega_\phi(u_1) \delta(v\gamma/u_2) +g_{q^{-1}}(v\gamma /u_2)\Omega_\psi (u_2)\delta(u_1\gamma/v)
 \\
  &\quad +g_{q^{-1}}(v\gamma /u_2)g(u_1\gamma/v)x^-(v)\Omega_\phi(u_1)\Omega_\psi (u_2)
\end{align*}
and
\begin{align*}
\Omega_\psi(u_2)\Omega_\phi(u_1)&x^-(v) =
  \Omega_\psi(u_2) \delta(u_1\gamma/v) +\Omega_\psi(u_2) x^-(v)g(u_1\gamma/v)\Omega_\phi(u_1)
  \\
  &= \Omega_\psi(u_2) \delta(u_1\gamma/v)  +g (u_1\gamma/v)\Omega_\phi(u_1)\delta(v\gamma/u_2)
\\
  &\quad+g_{q^{-1}}(v\gamma /u_2)g(u_1\gamma/v)x^-(v)\Omega_\psi(u_2)\Omega_\phi (u_1)
\end{align*}

Set $S=(q^2\gamma^2u_1-u_2)\Omega_\phi(u_1)\Omega_\psi(u_2)-(\gamma^2u_1-q^2u_2)\Omega_\psi(u_2)\Omega_\phi(u_1) $.
Then
\begin{align*}
Sx^-(v)
 &= (q^2\gamma^2u_1-u_2)\Omega_\phi(u_1) \delta(v\gamma/u_2)
    +(q^2\gamma^2u_1-u_2)g_{q^{-1}}(v\gamma /u_2)\Omega_\psi (u_2)\delta(u_1\gamma/v)
 \\
  &\quad +(q^2\gamma^2u_1- u_2)g_{q^{-1}}(v\gamma /u_2)g(u_1\gamma/v)x^-(v)\Omega_\phi(u_1)\Omega_\psi (u_2)  \\
  &\quad-( \gamma^2u_1-q^2u_2) \Omega_\psi(u_2) \delta(u_1\gamma/v)
    -( \gamma^2u_1-q^2u_2)g(u_1\gamma/v)\Omega_\phi(u_1)\delta(v\gamma/u_2)
\\
  &\quad-( \gamma^2u_1-q^2u_2)g_{q^{-1}}(v\gamma /u_2)g (u_1\gamma/v)x^-(v)\Omega_\psi(u_2)\Omega_\phi (u_1)\\  \\
&=\left((q^2\gamma^2u_1-u_2) -(\gamma^2u_1-q^2u_2)g (u_1\gamma/v)\right)\Omega_\phi(u_1)\delta(v\gamma/u_2)
 \\
  &\quad +\left((q^2\gamma^2u_1-u_2)g_{q^{-1}}(v\gamma /u_2)
    -(\gamma^2u_1-q^2u_2)\right) \Omega_\psi(u_2) \delta(u_1\gamma/v)
\\
  &\quad+g_{q^{-1}}(v\gamma /u_2)g (u_1\gamma/v)x^-(v)\\
  &\hskip 50pt \times\left((q^2\gamma^2u_1-u_2)\Omega_\phi(u_1)\Omega_\psi (u_2) -(\gamma^2u_1-q^2u_2)\Omega_\psi(u_2)\Omega_\phi (u_1)\right)\\  \\
  &=g_{q^{-1}}(v\gamma /u_2)g (u_1\gamma/v)x^-(v)S.
\end{align*}
As in the previous calculations we get that $S=0$ and thus the last statement of the proposition hold.

\end{proof}

The identities in \propref{commutatorprop} can be rewritten as
\begin{align}
(q^2v\gamma- u)\Omega_\psi(u)x^-(v)&=(q^2v\gamma- u)\delta(v\gamma/u)+(q^2v\gamma -u)x^-(v)\Omega_\psi(u),
\label{omegapsi2}\\
(q^2v- u\gamma)  \Omega_\phi(u)x^-(v)&=(q^2v- u\gamma)\delta(v/u\gamma)+( v-q^2u\gamma)x^-(v)\Omega_\phi(u)\label{omegaphi3}
\end{align}
which may be written out in terms of components as
\begin{align}
&q^2\gamma\Omega_\psi(m)x^-(n+1)- \Omega_\psi(m+1)x^-(n) \\
&\quad =(q^2\gamma-1)\delta_{m,-n-1}+ \gamma x^-(n+1)\Omega_\psi(m)-q^2x^-(n)\Omega_\psi(m+1),
\label{omegapsi4}\\
&  q^2\Omega_\phi(m)x^-(n+1)-  \gamma\Omega_\phi(m+1)x^-(n) \\
&\qquad =(q^2- \gamma)\delta_{m,-n-1}+ x^-(n+1)\Omega_\psi(m)-q^2\gamma x^-(n)\Omega_\psi(m+1),
\label{omegaphi5}
\end{align}

We can also write \eqnref{omegapsi} in terms of components and as operators on $\mathcal N_q^-$
\begin{equation}\label{omegapsi6}
    \Omega_\psi(k)x^-(m)=\delta_{k,-m}\gamma^{k}+\sum_{r\geq 0}g_{q^{-1}}(r)x^-(m+r)\Omega_\psi(k-r)\gamma^{r}.
\end{equation}
    The sum on the right hand side turns into a finite sum when applied to an element in $\mathcal N_q^-$, due to \eqnref{omegalocalfin}.

  We also have by \eqnref{omegaphipsi}
  \begin{equation}\label{omegaphipsi2}
    \Omega_\psi(k)\Omega_\phi(m)= \sum_{r\geq 0}g (r)\gamma^{2r}\Omega_\phi(r+m)\Omega_\psi(k-r),
\end{equation}
as operators on $\mathcal N_q^-$.

\section{The Kashiwara algebra $\mathcal K_q$}  The Kashiwara algebra $\mathcal K_q$ is defined to be the $\mathbb F(q^{1/2})$-algebra with generators $\Omega_\psi(m),x^-(n),\gamma^{\pm 1/2}$, $m,n\in\mathbb Z$ where $\gamma^{\pm 1/2}$ are central and the defining relations are
\begin{align}
q^2\gamma\Omega_\psi(m)&x^-(n+1)-  \Omega_\psi(m+1)x^-(n) \\
&=(q^2\gamma-1)\delta_{m,-n-1}+ \gamma x^-(n+1)\Omega_\psi(m)-q^2x^-(n)\Omega_\psi(m+1) \notag \\
q^2 \Omega_\psi(k+1)&\Omega_\psi(l) -
\Omega_\psi(l)\Omega_\psi(k+1)  =  \Omega_\psi(k)\Omega_\psi(l+1)
    - q^2\Omega_\psi(l+1)\Omega_\psi(k)\label{omegapsi3}
\end{align}
(which comes from \eqnref{omegapsi}, \eqnref{psipsi} written out in terms of components),
\begin{equation}
x^{-}(k+1)x^{-}(l) - q^{- 2}x^{-}(l)x^{-}(k+1)  = q^{- 2}x^{-}(k)x^{-}(l+1)
    - x^{-}(l+1)x^{-}(k)\label{xminusreln}
\end{equation}
 together with
\[
\gamma^{1/2}\gamma^{-1/2}=1=\gamma^{-1/2}\gamma^{1/2}.
\]

\begin{lem}  The $\mathbb F(q^{1/2})$-linear map $\bar{\alpha}:\mathcal K_q\to \mathcal K_q$ given by
$$
\bar{\alpha}(\gamma^{\pm 1/2})=\gamma^{\pm 1/2},\quad \bar{\alpha}(x^-(m))=\Omega_\psi (-m),\quad \bar{\alpha}(\Omega_\psi (m))=x^-(-m)
$$
for all $m\in\mathbb Z$ is an involutive anti-automorphism.
\end{lem}
\begin{proof}
We have
\begin{align*}
\bar{\alpha}&\left(x^{-}(k+1)x^{-}(l) - q^{- 2}x^{-}(l)x^{-}(k+1)  \right) \\
&=\Omega_\psi(-l)\Omega_\psi(-k-1) -q^{-2}\Omega_\psi(-k-1)\Omega_\psi(-l)  \\
&=q^{-2}\Omega_\psi(-l-1)\Omega_\psi(-k) -\Omega_\psi(-k)\Omega_\psi(-l-1)  \\
&=\bar{\alpha}\left(q^{-2}x^-(k)x^-(l+1)-x^-(l+1)x^-(k)\right)
\end{align*}
and
\begin{align*}
\bar{\alpha}&\left(q^2\gamma\Omega_\psi(m)x^-(n+1)-\Omega_\psi(m+1)x^-(n)\right) \\
    &=q^2\gamma\Omega_\psi(-n-1)x^-(-m)-\Omega_\psi(-n)x^-(-m-1)  \\
    &=(q^2\gamma-1)\delta_{-m,n+1}+ \gamma x^-(-m)\Omega_\psi(-n-1)-q^2x^-(-m-1)\Omega_\psi(-n) \\
    &=\bar{\alpha}\left((q^2\gamma-1)\delta_{m,-n-1}+ \gamma x^-(n+1)\Omega_\psi(m)-q^2x^-(n)\Omega_\psi(m+1)\right)
\end{align*}
\end{proof}

\begin{lem} $\mathcal N_q^-$ is a left $\mathcal K_q$-module.
\end{lem}
\begin{proof}  This follows from \eqnref{commutatorprop}
\end{proof}
\begin{lem}
$\mathcal N_q^-\cong \mathcal K_q/\sum_{k\in\mathbb Z}\mathcal K_q\Omega_\psi(k)$
\end{lem}
\begin{proof}
We have an induced left $\mathcal K_q$-module epimomorphism from $\mathcal K_q$ to $\mathcal N_q^-$ which sends $1$ to $1$.  Since $\Omega_\psi(k)$ annihilates $1$ for all $k$, we get an induced left  $\mathcal K_q^-$-module epimomorphism
\[
\begin{CD}
\mathcal K_q/\sum_k\mathcal K_q\Omega_\psi(k)@>\eta>> \mathcal N_q^-
\end{CD}
\]
Let  $C$ denote the subalgebra of $\mathcal K_q$ generated by $x^-(m),\gamma^{\pm1/2}$.  Then we have a surjective homomorphism
\[
\begin{CD}
C @>\mu>>\mathcal K_q/\sum_k\mathcal K_q\Omega_\psi(k)
\end{CD}
\]
The composition $\eta\circ \mu$ is surjective and since $\mathcal N_q^-$ is defined by generators $x^-(n),\gamma^{\pm 1/2}$ and relations \eqnref{Serre}, we get an induced map $\nu:\mathcal N_q^-\to C$ splitting the surjective map $\eta\circ \mu$. Since the composition $\nu\circ \eta\circ \mu$ is the identity, we get that $\eta\circ\mu$ is an isomorphism and thus $\eta$ is an isomorphism.
\end{proof}
\begin{prop}\label{form}
There is a unique symmetric form $(\enspace, \enspace)$ defined on $\mathcal N^-_q$ satisfying
$$
(x^-(m)a,b)=(a,\Omega_\psi(-m)b),\quad (1,1)=1.
$$
\end{prop}
\begin{proof}  Using the anti-automorphism $\bar\alpha$ we can make $M=\text{Hom}(\mathcal N^-_q,\mathbb F(q^{1/2}))$ into a left $\mathcal K_q$-module by defining
\begin{gather*}
(x^-(m)\phi)(a)=\phi(\Omega_\psi(-m)a),\enspace (\Omega_\psi(m)\phi)(a)=\phi(x^-(-m)a), \\
(\gamma^{\pm 1/2}\phi)(a)=\phi(\gamma^{\pm 1/2}a).
\end{gather*}
for $a\in \mathcal N_q^-$ and $\phi\in M$.

Consider the element $\beta_0\in M$ satisfying $\beta_0(1)=1$ and
$$
\beta_0\left(\sum_{k\in\mathbb Z}x^-(m)\mathcal K_q\right)=0.
$$
Then $\Omega_\psi(m)\beta_0=0$ for any $m\in\mathbb Z$, we get an induced homomorphism
$$
\bar\beta:\mathcal N_q^-\cong \mathcal K_q/\sum_{m\in\mathbb Z}\mathcal K_q\Omega_\psi(m)\to M.
$$
Define the bilinear form $(\enspace,\enspace):\mathcal K_q\times \mathcal K_q:\to \mathbb F(q^{1/2})$ by
$$
(a,b)=(\bar\beta(a))(b)
$$
This form satisfies $(1,1)=1$ and
\begin{gather*}
(x^-(m)a,b)=(a,\Omega_\psi(-m)b),\quad (\Omega_\psi(m)a,b)=(a,x^-(-m)b), \\ (\gamma^{\pm/2}a,b)=(a,\gamma^{\pm    /2}b).
\end{gather*}
Since $\mathcal N_q^-$ is generated by $x^-(m)$ and $\gamma^{\pm 1/2}$ we get that the form is the unique form satisfying these three conditions.  The form is symmetric since the form defined by $(a,b)'=(b,a)$ also satisfies the above conditions.
\end{proof}

\section{Imaginary Verma modules}
Let $\Lambda$ denotes the weight lattice of $A_1^{(1)}$,
$\lambda\in \Lambda$. Denote by $I^q(\lambda)$ the ideal of
$U_q=U_q(\hat{\mathfrak sl}(2))$ generated by $x^+(k)$, $k\in
\mathbb Z$, $a(l), l>0$, $K^{\pm 1}-q^{\lambda(h)}1$, $\gamma^{\pm
\frac{1}{2}}-q^{\pm \frac{1}{2}\lambda(c)}1$ and $D^{\pm 1}-q^{\pm
\lambda(d)}1$. The imaginary Verma module with highest weight
$\lambda$ is defined to be (\cite{MR97k:17014}) $$M_q(\lambda)=U/I^q(\lambda).$$

\begin{thm}[\cite{MR97k:17014}, Theorem 3.6] Imaginary Verma module $M_q(\lambda)$ is simple
 if and only if $\lambda(c)\neq 0$.
\end{thm}

Suppose now that $\lambda(c)=0$. Then $\gamma^{\pm \frac{1}{2}}$ acts on $M_q(\lambda)$ by $1$.  Consider an ideal
$J^q(\lambda)$ of  $U_q$ generated by $I^q(\lambda)$ and  $a(l)$ for
all $l$. Denote
$$\tilde{M}_q(\lambda)=U_q/J^q(\lambda).$$
Then $\tilde{M}_q(\lambda)$ is a homomorphic image of
 $M_q(\lambda)$ which we call \emph{reduced imaginary Verma
 module}. Module $\tilde{M}_q(\lambda)$ has a $\Lambda$-gradation:
 $$\tilde{M}_q(\lambda)=\sum_{\xi\in\Lambda}\tilde{M}_q(\lambda)_{\xi}.$$
 If $\alpha$ denotes a simple root of ${\mathfrak sl}(2)$ and $\delta$ denotes an indivisible imaginary root
 then $\tilde{M}_q(\lambda)_{\lambda-\xi}\neq 0$ if and only if $\xi=0$ or $\xi=-n\alpha+m\delta$ with $n>0$, $m\in \mathbb Z$.

 If $\xi=-n\alpha+m\delta$ then we set $|\xi|=n$. Note that $\mathcal N_q^-$ has also a $\Lambda$-grading:
 $x^-(n_1)x^-(n_2)\ldots x^-(n_k)\in (\mathcal N_q^-)_{\xi}$, where $\xi=-k\alpha+(n_1+\ldots+ n_k)\delta$, $|\xi|=k$.

In this section we discuss the properties of the reduced imaginary
Verma modules.

\begin{lem}\label{le-length2}
Let $\lambda\in \Lambda$ such that $\lambda(c)=\lambda(h)=0$,
$v\in\tilde{M}_q(\lambda)$ a nonzero element,
$v=u\tilde{v}_{\lambda}$, where $u\in (\mathcal N_q^-)_{\xi}$,
$|\xi|=2$. Then there exists $s\in \mathbb Z$ such that
$x^+(s)v\neq 0$.
\end{lem}

\begin{proof}
Let $\xi=-2\alpha+m\delta$, $m\in \mathbb Z$. We may assume
$$
u=\sum_{l}A_lx^-(l)x^-(m-l),
$$
where all but finitely many of the $A_l\in \mathbb C(q)$ are nonzero. Then by \eqnref{xcommutator}, we have
\begin{align*}
x^+(s)v=[x^+,u]\tilde{v}_{\lambda}=\frac{1}{q-q^{-1}}\sum_{l}A_l  \psi(s+l)x^-(m-l)\tilde{v}_{\lambda}
\end{align*}
for $s\gg 0$ as $\phi(l+s)=0$ for $l+s>0$ and $\psi(m-l+s)\tilde{v}_{\lambda} =0$ for $s\gg 0$.

Observe that
\begin{align*}
\psi(r)&=K(q-q^{-1})\Big(a(r)+\frac{1}{2}(q-q^{-1})\sum_{k_1+k_2=r\atop k_i\geq 1}
a(k_1)a(k_2) \\
&\qquad +\frac{1}{3!}(q-q^{-1})^2\sum_{k_1+k_2+k_3=r\atop k_i\geq 1}a(k_1)a(k_2)a(k_3)+ \ldots+\frac{1}{r!}
(q-q^{-1})^{r-1}a(1)^r)\Big),
\end{align*}
and for $k_1+\cdots +k_n=s+l$
$$
a(k_1)\cdots a(k_n)x^-(m-l)\tilde{v}_{\lambda} =(-1)^n \prod_{i=1}^n\frac{[2k_i]}{k_i}x^{-}(s+m)\tilde{v}_{\lambda} ,\
$$
by \eqnref{axcommutator}.
Thus we have
$$
x^+(s)v=\sum_l  A_l q^{2l+2s}f_l(s)Kx^-(s+m)\tilde{v}_{\lambda} ,
$$
and where for $s+l\geq 1$ one has
\begin{align*}
f_l(s)=
 \Big(&-\frac{[2(s+l)]}{s+l}+\frac{1}{2}t\sum_{k_1+k_2=s+l}
\frac{[2k_1][2k_2]}{k_1k_2} \\
&\quad -\frac{1}{3!}t^2\sum_{k_1+k_2+k_3=s+l}\frac{[2k_1][2k_2][2k_3]}{k_1k_2k_3}+\ldots +\frac{(-1)^{s+l}[2]^{s+l}}{(s+l)!}
t^{s+l-1}\Big),
\end{align*}
where
$t=q-q^{-1}$.
Note by \eqnref{quteformula} we have
 $$f_l(s)=\frac{(1-q^4)}{q^{2(s+l)}(q-q^{-1})}.$$

Suppose $x^+(s)v=0$ for any $s$. Then
$$
\sum_l  A_lq^{2l+2s}
f_l(s)=0
$$ for any sufficiently large $s$  and  so
\begin{equation}\label{confluencerln}
\sum_{l}A_l=0.
\end{equation}
  Note that this equality does not depend on $s$.

We can assume by \eqnref{Serre} that without loss of generality the monomials in $u$ are
ordered in such a way that $m-l\leq l$ for each $l$.
 Choose now the
smallest among $m-l$, say $r$, with $A_l=A_{m-r}\neq 0$ and apply $x^+(-r)$ to $\tilde{v}_\lambda$ noting that $r\leq l$ (so $l-r\geq 0$ and $-r+m-l\geq 0$):

\begin{align*}
x^+(-r)v&=\frac{1}{q-q^{-1}}\sum_{l}A_l
\psi(-r+l)x^-(m-l)\tilde{v}_{\lambda}+
A_{m-r}x^-(m-r)x^+(-r)x^-(r)\tilde{v}_{\lambda}\\
&=\sum_l  A_l
q^{2l-2r}f_l(-r)Kx^-(-r+m)\tilde{v}_{\lambda}+A_{m-r}x^-(m-r)\left(\frac{K-K^{-1}}{q-q^{-1}}\right)\tilde{v}_{\lambda} \\
&=A_{m-r}x^-(m-r)\left(\frac{K-K^{-1}}{q-q^{-1}}\right)\tilde{v}_{\lambda},
\end{align*}
due to \eqnref{confluencerln}.

This is a contradiction. It implies $v=0$.


\end{proof}

\begin{thm}\label{thm-imag-irred}
Let $\lambda\in \Lambda$ such that $\lambda(c)=0$. Then module $\tilde{M}_q(\lambda)$
is simple if and only if $\lambda(h)\neq 0$.
\end{thm}

\begin{proof}
Suppose $\lambda(h)=0.$ Let $v=x^-(m)\tilde{v}_{\lambda}$. Then for any $s\neq -m$ we have
$x^+(s)v=0$. Similarly,
$x^+(-m)v=[x^+(-m), x^-(m)]\tilde{v}_{\lambda}=\frac{1}{q-q^{-1}}(K-K^{-1})\tilde{v}_{\lambda}=0$, since
$K\tilde{v}_{\lambda}=q^{\lambda(h)}\tilde{v}_{\lambda}$. Hence, $v$ generates a proper nonzero submodule of
$\tilde{M}_q(\lambda)$.

Assume now $\lambda(h)\neq 0$.  To show simplicity of $\tilde{M}_q(\lambda)$ consider an arbitrary  homogeneous element
$v\in \tilde{M}_q(\lambda)_{\lambda-\xi}$ such that $x^+(s)v=0$ for any $s.$  We need to show that $v$ is a scalar multiple of $\tilde{v}_{\lambda}$.
We will proceed by the induction in $|\xi|$ to show that if $|\xi|>0$ then $v=0$.

 If  $|\xi|=1$ and $v\neq 0$ then $v=x^-(m)\tilde{v}_{\lambda}$ and
$x^+(-m)v\neq 0$. Hence $v=0$.  The case $|\xi|=2$ follows from
Lemma~\ref{le-length2}. Note that this case does not depend on the
value $\lambda(h)$.

Suppose now $|\xi|=k>2$, $v=u\tilde{v}_{\lambda}$ and
$$u=\sum_{n_1,\ldots, n_k}A(n_1, \ldots, n_k)x^-(n_1)\ldots x^-(n_k).$$
Using notation from the lemma above we have
$$
\psi(s)x^-(n)=-tf_n(s)Kx^-(n+s)+q^{-2}x^-(n)\psi(s)
$$

and thus
\begin{align*}
x^+(s)&(x^-(n_1)\cdots x^-(n_k)\tilde{v}_{\lambda})=\frac{1}{q-q^{-1}}\psi(n_1+s)x^-(n_2)\cdots x^-(n_k)\tilde{v}_{\lambda} \\
&\quad +x^-(n_1)x^+(s)x^-(n_2)\cdots x^-(n_k)\tilde{v}_{\lambda} \\
&=-q^{2s+2n_1}f_{n_1}(s+n_1)Kx^-(n_1+n_2+s)x^-(n_3)\cdots x^-(n_k)\tilde{v}_{\lambda} \\
&\quad +\frac{q^{-2}}{q-q^{-1}}x^-(n_2)\psi(n_1+s)x^-(n_3)\cdots x^-(n_k)\tilde{v}_{\lambda} \\
&\quad +x(n_1)x^+(s)x^-(n_2)\cdots x^-(n_k)\tilde{v}_{\lambda} \\ \\
&=-q^{2s+2n_1}f_{n_1}(s)Kx^-(n_1+n_2+s)x^-(n_3)\cdots x^-(n_k)\tilde{v}_{\lambda} \\
&\quad -q^{-2}q^{2s+2n_1}f_{n_1}(n_1+s)x^-(n_2)Kx^-(n_1+n_3+s)x^-(n_4)\cdots x^-(n_k)\tilde{v}_{\lambda}\\
&\quad +\frac{q^{-4}}{q-q^{-1}}x^-(n_2)x^-(n_3)\psi(n_1+s)x^-(n_4)\cdots
x^-(n_k)\tilde{v}_{\lambda} \\
&\quad+\frac{1}{q-q^{-1}}x^-(n_1)\psi(n_2+s)x^-(n_3)\cdots x^-(n_k) \\
&\quad + x^-(n_1)x^-(n_2)x^+(s)x^-(n_3)\cdots x^-(n_k)\tilde{v}_{\lambda} \\
&\quad-q^{2s+2n_{k-1}}f_{n_{k-1}}(s+n_{k-1})x^-(n_1)\cdots x^-(n_{k-2})Kx^-(n_{k-1}+n_k+s)\tilde{v}_{\lambda}.
\end{align*}

We may order monomials in $u$ in such a way that $n_1\leq n_2\leq
\ldots \leq n_k$. We also introduce  lexicographical ordering
among the monomials.

 The smallest monomial in the image $x^+(s)(x(n_1)\ldots
 x^-(n_k)\tilde{v}_{\lambda})$ is
 $$
 x^-(n_1)\ldots
 x^-(n_{k-2})Kx^-(n_{k-1}+n_k+s)
 $$
 up to a constant. It determines
 uniquely the first $k-2$ elements in the monomial and leaves a
 freedom in the choice of last two elements (remembering that $u$ is homogeneous).
 Hence, we may assume that
 $$
 u=x^-(n_1)x^-(n_2)\ldots x^-(n_{k-2})\sum_l B_l x^-(m-l)x^-(l),
 $$
 for some fixed $m$, $n_1\leq \cdots\leq n_{k-2}\leq m-l\leq l$. Then
\begin{align*}
x^+(s)v&=x^+(s)u\tilde{v}_{\lambda}  \\
  &=[x^+(s), x^-(n_1)x^-(n_2)\ldots x^-(n_{k-2})]\sum_l B_lx^-(m-l) x^-(l)\tilde{v}_{\lambda} \\
&\quad +x^-(n_1)x^-(n_2)\ldots
  x^-(n_{k-2})[x^+(s),\sum_l B_lx^-(m-l) x^-(l)]\tilde{v}_{\lambda}.
  \end{align*}
  Note that the first part in the sum above will contribute
  smaller monomials than the second part.   Hence, if $x^+(s)v=0$
  for any $s\in \mathbb Z$ then $$
  [x^+(s),\sum_l B_l x^-(l)x^-(m-l)]\tilde{v}_{\lambda}=0,
  $$
  for all sufficiently large integers $s$.  Define $A_l\in\mathbb C(q)$ such that
  $$
  \sum_l A_lx^-(l)x^-(m-l) =\sum_l B_lx^-(m-l) x^-(l)
  $$

Applying Lemma~\ref{le-length2}  we obtain that all $A_l$ are zero (and hence so are the $B_l$)
and thus $v=0$. This completes the proof.
\color{black}
\end{proof}

 Set
$R_q(\lambda)=\sum_{\xi, |\xi|>0}\tilde{M}_q(\lambda)_{\xi}$. Then
$R_q(\lambda)$ is the unique maximal submodule of
$\tilde{M}_q(\lambda)$ and $\dim
\tilde{M}_q(\lambda)/R_q(\lambda)=1$.




\begin{rem}
It was shown in \cite{MR97k:17014}, Theorem 5.4 that imaginary Verma
module $M(\lambda)$ over affine $\hat{\mathfrak sl}(2)$ admits a
quantum deformation to the imaginary Verma module $M_q(\lambda)$
over $U_q$ in such a way that the dimensions of the weight spaces
are invariant under the deformation, generalizing the Lusztig's
deformation functor constructed originally for classical Verma
modules \cite{MR954661}, see also \cite{MR1662112}.
Theorem~\ref{thm-imag-irred} shows that Lusztig's deformation
functor  preserves module structure in the case of imaginary Verma
modules (see \cite{MR95a:17030}).
\end{rem}

\section{Simplicity of $\mathcal N_q^-$ as a $\mathcal K_q$-module}

We will show that $\mathcal N_q^-$ is simple as a module over
$\mathcal K_q$.

\begin{lem}\label{le-primitive}
Let $P\in \mathcal N_q^-$. If $\Omega_{\psi}(s)P=0$ for any $s\in
\mathbb Z$, then $P$ is a constant multiple of $1$.

\end{lem}

\begin{proof}
We may assume without loss of generality that $P$ is a homogeneous
element, say $P\in (\mathcal N_q^-)_{\lambda-\xi}$. We assume
that $\xi\neq 0$. Then $\xi=n\alpha+m\delta$, $n>0$, $m\in \mathbb
Z$. Set $|\xi|=n$. We shall prove the lemma by  induction on
$|\xi|$.

Suppose $|\xi|=1$. Then $P=x^-(m)$ and
\begin{align*}
\Omega_{\psi}(s)(P)&=\delta_{s,-m}\gamma^s+\sum_{r\in \mathbb Z}g(s-r)x^-(m-r+s)\gamma^r\Omega_{\psi}(r)1 \\
&=\delta_{s,-m}.
\end{align*}
Hence $\Omega_{\psi}(-m)(P)\neq 0$ unless $P=0$.

Suppose $|\xi|>1$.  We assume $\Omega_\psi(l)(P)=0$ for any $l\in\mathbb Z$ and then we use \eqnref{omegaphipsi2}.  For all $k$ and $m$ we get

  \begin{equation}
    \Omega_\psi(k)\Omega_\phi(m)(P)= \sum_{r\geq 0}g (r)\gamma^{2r}\Omega_\phi(r+m)\Omega_\psi(k-r)(P)=0,
\end{equation}
Hence by the induction hypothesis $\Omega_\phi(m)(P)=0$ as $\Omega_\phi(m)(P) \in (\mathcal N_q^-)_{\lambda-\xi+1}$.
Then $[x^+(m),P]=0$ by \eqnref{xplusP}.

Consider the imaginary Verma module
$M_q(\lambda)$ with $\lambda(c)=0$ and choose $\lambda$ such that $\lambda(h)\neq 0$. Then $\tilde{M}_q(\lambda)$ is the unique irreducible quotient of $M_q(\lambda)$
 and $v=Pv_{\lambda}$ is a nonzero element of the module $\tilde{M}_q(\lambda)$.

\color{black}
Thus
$$x^+(s)v=[x^+, P]\tilde{v}_{\lambda}+ Px^+(s) \tilde{v}_{\lambda}=0$$
for all $s\in \mathbb Z$.

Consider $V= \mathcal N_q^-v\subset \tilde{M}_q(\lambda)$. Then $V$ is a nonzero proper submodule
of $\tilde{M}_q(\lambda)$ which is a contradiction by Theorem~\ref{thm-imag-irred}.
This completes the proof.
\end{proof}

\begin{rem}
Suppose $|\xi|=2$. We will give a direct proof of
Lemma~\ref{le-primitive} in this case without the use of
Theorem~\ref{thm-imag-irred}.

Let
$$
P=\sum_{n_1, n_2,
n_1+n_2=m}A(n_1,n_2)x^-(n_1)x^-(n_2).
$$
We can assume that $n_1\leq
n_2$ in all the monomials in $P$. Then
\begin{align*}
\Omega_{\psi}(s)(P)&=\sum_{n_1, n_2,
n_1+n_2=m}A(n_1,n_2)(\delta_{n_1,-s}x^-(n_2)+g(n_2+s)x^-(m+s))\gamma^s  \\
&=A(-s, m+s)\delta_{n_1,-s}x^-(m+s)\gamma^s \\
&\hskip 50pt+\left(\sum_{n_1}A(n_1, m-n_1)g(m+s-n_1)\right)x^-(m+s)\gamma^s.
\end{align*}
If $\Omega_{\psi}(s)(P)=0$ for all $s$ then, in particular,
$$\sum_{n_1}A(n_1, m-n_1)g(m+s-n_1)=0$$
for any $s$ sufficiently large. Since $g(p)=(q^4-1)q^{-2p-2}$ for
$p\geq 1$, we will get

$$(q^4-1)q^{-2(m+s+1)}\sum_{n_1}q^{2n_1}A(n_1, m-n_1)=0$$
implying $\sum_{n_1}q^{2n_1}A(n_1, m-n_1)=0$.  Note that this
relation does not depend on $s$.

Choose now the smallest among $n_1$, say $n_1=r$, with $A(n_1,m-n_1)\neq 0$ and apply
$\Omega_{\psi}(-r)$:
$$
\Omega_{\psi}(-r)P=
A(r,m-r)x(m-2r) \gamma^s+ \sum_{n_1}A(n_1,
m-n_1)g(m-n_1-r)x^-(m-2r)\gamma^s=0,
$$
Now
$$
 \sum_{n_1}A(n_1,m-n_1)g(m-n_1-r)=(q^4-1)q^{-2m+2r-2} \sum_{n_1}A(n_1,m-n_1)q^{2n_1}=0.
$$
Hence
$$
0=\Omega_{\psi}(-r)P=
A(r,m-r)x(m-2r) \gamma^s,
$$
but then $A(r,m-r)=0$ which is a contradiction.

We have
Suppose $m-2r>0$. Then
 $$\sum_{n_1}A(n_1, m-n_1)g(m-n_1-r)=(q^4-1)q^{-2(m-r+1)}\sum_{n_1}q^{2n_1}A(n_1,m-n_1)=0,$$

and
$$\sum_{n_1}A(n_1, m-n_1)g(m-n_1-r)x^-(m-2r)\gamma^s=0.$$

Thus $A(r, m-r)=0$, which is a contradiction.

If $m=2r$ then we have a unique monomial $x^-(r)^2$ in $P$ due to the chosen ordering. Hence,
$$\Omega_{\psi}(-r)P=
A(r,r)(1+g(0))x^-(0) \gamma^s=A(r,r)(1+ q^{-2})x^-(0) \gamma^s=0$$
implies $A(r,r)=0$.
This is again a contradiction.
 Therefore, there exists $s$ such that
$\Omega_{\psi}(s)P\neq 0$.  Note that in fact we proved that in
the case $|\xi|=2$, $\Omega_{\psi}(s)P\neq 0$ for all $s\in
\mathbb Z$.

\end{rem}

 Lemma~\ref{le-primitive} implies immediately the following result.

\begin{thm}
The algebra $\mathcal N_q^-$ is simple as a $\mathcal K_q$-module.
\end{thm}

\begin{cor}  The form $(\enspace,\enspace)$ defined in \propref{form} is non-degenerate.
\end{cor}
\begin{proof}  By \propref{form} the radical of the form $(\enspace,\enspace)$ is a  $\mathcal K_q$-submodule of $\mathcal N_q^-$ and since $(1,1)=1$, the radical must be zero.
\end{proof}

\section{Acknowledgement}
The authors are grateful to the organizers for the invitation to the conference at Banff where this project was initiated. The first author would like to thank North Carolina State University for the support and hospitality during his numerous visits to Raleigh.  The second author was partially supported by Fapesp (processo 2005/60337-2) and CNPq (processo 301743/2007-0). He is grateful to the North Carolina State University for the support and hospitality during his visit to Raleigh. The third author was partially supported by the NSA grant H98230-08-1-0080.


\def\cprime{$'$}
\providecommand{\bysame}{\leavevmode\hbox to3em{\hrulefill}\thinspace}
\providecommand{\MR}{\relax\ifhmode\unskip\space\fi MR }
\providecommand{\MRhref}[2]{%
  \href{http://www.ams.org/mathscinet-getitem?mr=#1}{#2}
}
\providecommand{\href}[2]{#2}

\end{document}